\documentclass[12pt]{article}
\usepackage{amssymb,a4}
\usepackage{amsmath,amsfonts,amssymb,amsthm}
\def\Res{{\rm Res}}
\def\wt{{\rm wt}}

\def\C{{\mathbb C}}

\def\Z{{\mathbb Z}}

\def\N{{\mathbb N}}
\def\1{{\bf 1}}
\def\Hom{{\rm Hom}}
\def\End{{\rm End}}

\def\theequation{5.\arabic{equation}}

\renewcommand{\theequation}{\thesection.\arabic{equation}}

\newtheorem{theorem}{Theorem}[section]

\newtheorem{coro}[theorem]{Corollary}

\newtheorem{conjecture}[theorem]{Conjecture}
\theoremstyle{definition}

\begin{document}
\begin{center}
{\Large {\bf Representations of vertex operator algebras}} \\

\vspace{0.5cm} Chongying Dong\footnote{Supported by NSF grants,
China NSF grant 10328102 and  a Faculty research grant from  the
University of California at Santa Cruz.}
\\
Department of Mathematics\\ University of
California\\ Santa Cruz, CA 95064 \\
Cuipo Jiang\footnote{Supported in part by China NSF grant
10571119.}\\
 Department of Mathematics\\ Shanghai Jiaotong University\\
Shanghai 200030 China\\
\hspace{1cm}

This paper is dedicated to Professors James Lepowsky and Robert
Wilson on the occasion of their 60th birthday
\end{center}

\hspace{1cm}
\begin{abstract} This paper is an exposition of the representation
theory of vertex operator algebras in terms of associative
algebras $A_n(V)$ and their bimodules. A new result on the
rationality is given. That is, a simple vertex operator algebra
$V$ is rational if and only if its Zhu algebra $A(V)$ is a
semisimple associative algebra and each irreducible admissible
$V$-module is ordinary.
 2000MSC:17B69
\end{abstract}

\section{Introduction}
\def\theequation{1. \arabic{equation}}
\setcounter{equation}{0}

Rational vertex operator algebras, which play a fundamental role in
rational conformal field theory (see [BPZ] and [MS]), form an
important class of vertex operator algebras.
Most vertex operator algebras which have been studied so
far are rational vertex operator algebras. Familiar examples include
the moonshine module $V^{\natural}$ ([B], [FLM], [D2]), the vertex
operator algebras $V_L$ associated with positive definite even
lattices $L$ ([B], [FLM], [D1]), the vertex operator algebras
associated with integrable representations of affine Lie algebras
([DL], [FZ]), the vertex operator algebras $L(c_{p,q},0)$
associated with irreducible highest weight representations for the
discrete series of the Virasoro algebra ([DMZ], [W1]), and framed
vertex operator algebras ([DGH], [M]).

The notion of rational vertex operator algebra was first
introduced in \cite{Z} to study the modular invariance of trace
functions for vertex operator algebras. The rational vertex
operator algebra defined in \cite{Z} satisfies the following three
conditions: (1) Any admissible module is completely reducible. (2)
There are only finitely many irreducible admissible modules up to
isomorphisms. (3) Each irreducible admissible module is ordinary.
It was proved in \cite{DLM2} that conditions (2) and (3) are
consequences of condition (1). This makes the notion of
rationality more natural comparing to the notion of semisimplicity
of a finite dimensional associative algebra.

The rationality given here is clearly an external definition. One
needs to verify the complete reducibility for any admissible
module. So in practice, it is very hard to check if a vertex
operator algebra is rational. In the classical case for finite
dimensional associative algebra, there is an internal
characterization of semisimplicity. That is, the Jacobson or
solvable radical of the algebra is zero. Similarly, for a finite
dimensional Lie algebra, one can use the Killing form or solvable
radical to define semisimple Lie algebras. A natural question is
whether there is an internal definition of rationality for vertex
operator algebras.

There is a close relation between the representation theory of a
vertex operator algebra $V$ and the representation theory of
$A(V)$ which is an associative algebra associated to $V$ defined
in \cite{Z}. The most important result about $A(V)$ is that there
is a one to one correspondence between the equivalence classes of irreducible
admissible $V$-modules and the equivalence classes of simple $A(V)$-modules.
This turns $A(V)$ into a very powerful and effective tool in
the classification of irreducible modules for a given vertex operator
algebra. In fact, the irreducible modules for many families of vertex
operator algebras such as affine vertex operator algebras
\cite{FZ}, the Virasoro vertex operator algebras \cite{DMZ}, \cite{W1},
lattice type vertex operator algebras \cite{DN1}-\cite{DN3},
\cite{AD}, \cite{TY}, certain $W$-algebras \cite{W2}, \cite{DLTYY}
and some other vertex operator algebras \cite{A2}, \cite{Ad},
\cite{KW}, \cite{KMY}  have been classified using the
$A(V)$-theory.

Motivated by the representation theory of vertex operator algebras,
a sequence of associative algebras $A_n(V)$ for $n\geq 0$ were
introduced in \cite{DLM3} to deal with  the first $n+1$
homogeneous subspaces of an admissible module such that $A_0(V)$
is exactly $A(V).$ It is proved in \cite{DLM3} that $V$ is
rational if and only if $A_n(V)$ are finite dimensional semisimple
associative algebras for all $n.$ Since the construction of
$A_n(V)$ is very complicated, it is hard to compute $A_n(V)$ for
large $n.$ It has been suspected for a long time
 that the semisimplicity of $A(V)$ is good enough
to characterize the rationality of $V.$

This paper is an exposition and extension of our recent results in
\cite{DJ} concerning the rationality of a vertex operator algebra
together with the $A_n(V)$ theory. We obtain that if $V$ is a
simple vertex operator algebra such that each weight space of any
irreducible admissible $V$-module is finite dimensional, then $V$
is rational if and only if $A(V)$ is semisimple. In particular, if
$V$ is $C_2$-cofinite then the semisimplicity of $A(V),$
rationality of $V$ and regularity of $V$ are all equivalent. This
surely gives an internal characterization of rationality as the
semisimplcity of $A(V)$ can be defined internally. The main tool
used to prove this result is the bimodule theory developed in
\cite{DJ}.

Since the weight spaces of an irreducible admissible $V$-module
for all well known vertex operator algebras are
finite-dimensional,  the assumption on the dimensions of weight
spaces is not a strong assumption. For most vertex operator algebras,
the computation of $A(V)$ is necessary to classify the
irreducible modules. So in practice, this new result is very
useful. For example, the rationality of $V_L^+$ for a rank one
positive definite lattice $L$ proved in \cite{A1} is an easy
corollary of our result.

We firmly believe that a simple vertex operator algebra is rational
if and only if $A(V)$ is semisimple without any assumption. But
we cannot achieve this in this paper.

\section{Preliminaries}
\def\theequation{2.\arabic{equation}}
\setcounter{equation}{0}

Let $V=(V,Y,1,\omega)$ be a vertex operator algebra (see \cite{B} and
\cite{FLM}). We define weak, admissible and ordinary modules following
\cite{DLM1}-\cite{DLM2}. We also define rationality, $C_2$-cofiniteness
and regularity.

 A {\em weak module}  $M$ for $V$ is a vector space equipped with a linear map
$$\begin{array}{l}
V\to (\mbox{End}\,M)[[z^{-1},z]]\label{map}\\
v\mapsto\displaystyle{Y_M(v,z)=\sum_{n\in \Z}v_nz^{-n-1}\ \ (v_n\in\mbox{End}\,M)}
\mbox{ for }v\in V\label{1/2}
\end{array}$$
(where for any vector space $W,$ we define $W[[z^{-1},z]]$ to be the vector
space of $W$-valued formal series in $z$)
satisfying the following conditions for $u,v\in V$,
$w\in M$:
\begin{eqnarray*}\label{e2.1}
& &v_nw=0\ \ \
\mbox{for}\ \ \ n\in \Z \ \ \mbox{sufficiently\ large};\label{vlw0}\\
& &Y_M({\bf 1},z)=1;\label{vacuum}
\end{eqnarray*}
$$\begin{array}{c}
\displaystyle{z^{-1}_0\delta\left(\frac{z_1-z_2}{z_0}\right)
Y_M(u,z_1)Y_M(v,z_2)-z^{-1}_0\delta\left(\frac{z_2-z_1}{-z_0}\right)
Y_M(v,z_2)Y_M(u,z_1)}\\
\displaystyle{=z_2^{-1}\delta\left(\frac{z_1-z_0}{z_2}\right)
Y_M(Y(u,z_0)v,z_2)}.
\end{array}$$
This completes the definition. We denote this module by $(M,Y_M)$
(or briefly by $M$). It is proved in \cite{DLM1} that the
following are true on the weak module $M:$
$$Y_M(L(-1)u,z)=\frac{d}{dz}Y_M(u,z)$$
and
$$[L(m),L(n)]=(m-n)L(m+n)+\frac{m^3-m}{12}\delta_{m+n,0}C$$
where $Y_M(\omega,z)=\sum_{n\in\Z}L(n)z^{-n-2}.$

An {\em ordinary} $V$-module is a $\C$-graded  weak $V$-module
$$M=\bigoplus_{\lambda \in{\C}}M_{\lambda} $$
such that $\dim M_{\lambda}$ is finite and $M_{\lambda+n}=0$ for
fixed $\lambda$ and $n\in {\Z}$ small enough where $M_{\lambda}$
is the $\lambda$-eigenspace for $L(0)$ with eigenvalue $\lambda.$

 An {\em admissible} $V$-module is
a  weak $V$-module $M$ which carries a
${\Z}_{+}$-grading
$$M=\bigoplus_{n\in {\Z}_{+}}M(n)$$
($\Z_+$ is the set of all nonnegative integers) such that if $r, m\in
{\Z} ,n\in {\Z}_{+}$ and $a\in V_{r}$ then $a_{m}M(n)\subseteq
M(r+n-m-1).$ Since the uniform degree shift gives an isomorphic
admissible module we assume $M(0)\ne 0$ in many occasions. The
admissible module is also called the $\N$-graded module in
\cite{LL}. It is easy to prove that any ordinary module is an
admissible module.

Following \cite{DLM2}, a vertex operator algebra is called {\em
rational} if any admissible module is a direct sum of irreducible
admissible modules. As we mentioned in the introduction,  the
rationality given here looks weaker than the rationality defined
originally in \cite{Z} although it was proved in \cite{DLM2} these
two rationalities are the same, also see Theorem  \ref{P3.1}
below.

Another important concept is the $C_2$-cofiniteness.
The $C_2$-cofiniteness is a very important concept in the theory
of vertex operator algebra. Many important results such as modular
invariance of trace functions [Z], [DLM5], and that the weight one
subspace of a CFT type vertex operator algebra is reductive \cite{DM1},
tensor product and Verlinde formula \cite{H}
 need both rationality and $C_2$-cofininteness. The $C_2$-cofinite
condition also plays a fundamental role in the study of conformal field
theory [GN], [NT]. In particular,
the $C_2$-cofiniteness implies that the fusion rules are finite
[Bu] and the conformal blocks are finite dimensional [NT].

A vertex operator algebra $V$ is called $C_2$-cofinite if the
subspace $C_2(V)$ spanned by $u_{-2}v$ for all $u,v\in V$ has
finite codimension in $V$ \cite{Z}. It is demonstrated in
\cite{GN} that if $V$ is $C_2$-cofinite then $V$ is finitely
generated with a PBW-like spanning set.
\begin{theorem}\label{GN} Let $X=\{x^\alpha|\alpha \in I\}$ be a
subset of $V$ consisting of homgeneous vectors such that
$\{x^\alpha+C_2(V)|\alpha \in I\}$ form a basis of $V/C_2(V).$
Then $V$ is spanned by
$$x^{1}_{-n_1} x^{2}_{-n_2} \cdots x^{k}_{-n_k} {\textbf 1}$$
where $n_1>n_2> \cdots >n_k > 0$ and $x^{i} \in X$ for $1 \leq i
\leq k$. In particular, $V$ is finitely generated if $V$ is
$C_2$-cofinite.
\end{theorem}

This result has been generalized in \cite{Bu} to get a PBW-like
spanning set for any weak module generated by one element. We
should remark that there are some other results in \cite{KL}
concerning finitely generated properties.

A vertex operator algebra is called {\em regular} if any weak
module is a direct sum of irreducible ordinary modules. Clearly, a
regular vertex operator algebra is rational. It has been a very
challenging problem to find the relationship among rationality,
regularity and $C_2$-cofiniteness. It is shown in \cite{ABD} and
\cite{L} that rationality together with $C_2$-cofiniteness and
regularity are equivalent. Although it has been suspected that
rationality and $C_2$-cofiniteness are equivalent, there is
a counter example in \cite{A2} where a vertex operator algebra is
$C_2$-cofinite but not rational. On the other hand, the $C_2$-cofiniteness
is good enough to guarantee the integrability of the vertex operator algebra
of CFT type \cite{DM2}. That is, the vertex operator subalgebra
of $V$ generated by any semisimple subalgebra of $V_1$ is integrable
and $V$ is an integrable module for the corresponding affine Kac-Moody
Lie algebra.

In the definition of regular vertex operator algebra, we require that any weak
module is a direct sum of irreducible ordinary modules. An immediate question
is the following: {\em If every weak module is a direct sum of irreducible
weak modules, is $V$ regular?} It was proved in \cite{ABD} that
if $V$ is $C_2$-cofinite then any finitely generated weak module
is ordinary. So the answer to this question would be positive
if one could prove that rationality implies $C_2$-cofiniteness. But this
again is a very difficult problem.

\section{The associative algebra $A(V)$}
\setcounter{equation}{0}

In this section we review the $A(V)$ theory following \cite{DLM2}
and \cite{Z}. Let $V$ be a vertex operator algebra. For any homogeneous vectors $a\in V$, and $b\in V$, we define
\begin{align*}
&a*b=\left(\Res_{z}\frac{(1+z)^{\wt{a}}}{z}
Y(a,z)\right)b,\\
&a\circ b=\left( \Res_{z}\frac{(1+z)^{\wt{a}}}{z^{2}}
Y(a,z)\right)b,
\end{align*}
and extend to $V$ bilinearly.  Denote by $O(V)$ the linear
span of $a\circ b$ ($a,b\in V$) and set $A(V)=V/O(V)$.

The definition of $A(V)$ is very natural from the
representation theory of vertex operator algebra. Let $M$ be an
admissible $V$-module. For each homogeneous $v\in V,$ we set
$o(v)=v_{\wt v-1}$  on $M$ and extend linearly to whole $V.$ Then
it is not hard to prove that $o(u)o(v)=o(u*v)$ on $M(0)$ for all
$u,v\in V.$ So the operation $*$ will be the product in $A(V).$
Since $o((L(-1)+L(0))v)=0$ on $M$ we have $o(((L(-1)+L(0))v)*u)=0$
on $M(0)$ for all $u,v\in V.$ So we have to modulo out all
$((L(-1)+L(0))v)*u$ in order to get an effective action. One can
verify that $((L(-1)+L(0))v)*u=v\circ u.$ This should explain why
the algebra $A(V)$ is so defined.

We write $[a]$ for the image $a+O(V)$ of $a\in V$.  The following theorem is due to
\cite{Z} and \cite{DLM2}.
\begin{theorem}\label{P3.1} Assume that $V$ is a vertex operator
algebra. Then

(1) The bilinear operation $*$ on $V$ induces an associative algebra
structure on $A(V).$ The vector $[\1]$ is the identity and $[\omega]$ is in the center of
$A(V)$.

(2) Let $M=\bigoplus_{n=0}^{\infty}M(n)$ be an admissible $V$-module with
$M(0)\ne 0.$
Then the linear map
\[
o:V\rightarrow\End M(0),\;a\mapsto o(a)|_{M(0)}
\]
induces an algebra homomorphism from $A(V)$ to $\End M(0)$.
Thus $M(0)$ is a left $A(V)$-module.

(3) The map $M\mapsto M(0)$ induces a bijection from the set of equivalence
classes of irreducible admissible $V$-modules to the set of equivalence classes of irreducible $A(V)$-modules.

(4) If $V$ is rational, then $A(V)$ is a finite dimensional semisimple
associative algebra.

(5) If $V$ is rational then there are only finitely many
irreducible admissible $V$-modules up to isomorphism and each
irreducible admissible module is ordinary.
\end{theorem}

We remark that Theorem \ref{P3.1} (5) was part of the definition of
rational vertex operator algebra given in \cite{Z}. If $V$ is
$C_2$-cofinite, then Theorem \ref{P3.1} (5) also holds (see
\cite{KL}). This suggests that there is a strong link between
rationality and $C_2$-cofiniteness. But it is not clear what the
link might be.

Theorem \ref{P3.1} is very powerful in the classification of
irreducible modules for a vertex operator algebra. In order to
classify the irreducible modules for a vertex operator algebra
$V$, it is enough to classify the irreducible modules for $A(V)$
which is computable in many cases. The classifications of
irreducible modules have been achieved for many well known vertex
operator algebras (see \cite{FZ}, \cite{W1}, \cite{W2}, \cite{KW},
\cite{DN1}-\cite{DN3}, \cite{Ad}, \cite{AD}, \cite{DLTYY},
\cite{TY}, \cite{KMY}, \cite{A2}) in this way. But the structure
of $A(V)$ cannot tell whether or not a vertex operator algebra $V$
is rational at this stage.

\section{The associative algebras $A_n(V)$}
\setcounter{equation}{0}

Although $A(V)$ theory is very useful, the $A(V)$ is not good
enough to determine the representation theory for $V$ completely.
Motivated by the ``graded structure'' of admissible modules,
 an associative algebra $A_n(V)$ was introduced in \cite{DLM3}
for any nonnegative integer $n,$ generalizing the Zhu algebra
$A(V)$ which is $A_0(V).$ The role of $A_n(V)$ playing in the
representation theory of vertex operator algebra will become clear
after Theorem \ref{tha} below.

Let $O_n(V)$ be the linear span of all $u\circ_n v$ and $L(-1)u+L(0)u$
where for homogeneous $u\in V$ and $v\in V,$
$$u\circ_n v=\Res_{z}Y(u,z)v\frac{(1+z)^{\wt u+n}}{z^{2n+2}}.
$$
Define the linear space $A_n(V)$ to be the quotient $V/O_{n}(V).$
We also define a second product $*_n$ on $V$ for $u$ and $v$ as
follows:
$$u*_nv=\sum_{m=0}^{n}(-1)^m{m+n\choose n}\Res_zY(u,z)\frac{(1+z)^{\wt\,u+n}}{z^{n+m+1}}v.$$
Extend linearly to obtain a bilinear product  on $V$ which
coincides with that of Zhu algebra [Z] if $n=0.$ We denote the
product $*_0$  by $*$ in this case.

Here are the main results on $A_n(V)$ obtained in \cite{DLM3}.
\begin{theorem}\label{tha} Let $V$ be a vertex operator algebra and $n$ a nonnegative
integer. Then

(1) $A_n(V)$ is an associative algebra whose product is induced by
$*_n.$

(2) The identity map on $V$ induces an algebra epimorphism from
$A_n(V)$ to $A_{n-1}(V).$

(3) Let $W$ be a weak module and set
$$\Omega_n(W)=\{w\in W|u_mw=0, u\in V, m> \wt u-1+n\}.$$
Then $\Omega_n(W)$ is an $A_n(V)$-module such that $v+O_n(V)$ acts as $o(v).$

(4) Let $M=\sum_{m=0}^{\infty}M(m)$ be an admissible $V$-module.
Then each $M(m)$ for $m\leq n$ is
an $A_n(V)$-submodule of $\Omega_n(W).$ Furthermore,
$M$ is irreducible if and only if each $M(n)$ is an irreducible
$A_n(V)$-module.

(5) For any $A_n(V)$-module $U$ which cannot factor through
$A_{n-1}(V)$ there is a unique Verma type admissible $V$-module
$M(U)$ generated by $U$ so that $M(U)(0)\ne 0$ and $M(U)(n)=U.$
Moreover, for any weak $V$-module $W$ and any $A_n(V)$-module
homomorphism $f$ from $U$ to $\Omega_n(W)$ there is a unique
$V$-homomorphism from $M(U)$ to $W$ which extends  $f.$

(6) $V$ is rational if and only if $A_n(V)$ are finite dimensional
semisimple algebras for all $n\geq 0.$
\end{theorem}

Theorem \ref{tha} is a generalization and extension of Theorem
\ref{P3.1}. In (5) the assumption that $U$ cannot factor through
$A_{n-1}(V)$ is not essential. In fact, for any $A_n(V)$-module
$U$ we can construct $M(U)$ but we cannot conclude that
$M(U)(0)\ne 0.$ The algebra $A_n(V)$ definitely carries more
information than $A(V)$ on the admissible $V$-modules. For
example, one cannot tell if an admissible $V$-module
$M=\sum_{n\geq 0}M(n)$ is irreducible even though $M(0)$ is an
irreducible $A(V)$-module. The $A_n(V)$ theory plays a fundamental
role in the study of orbifold theory and dual pairs (see
\cite{DLM0}, \cite{DY}, \cite{Y}, \cite{MT}).

\section{Bimodules $A_{n,m}(V)$}
\setcounter{equation}{0}

Let $M=\sum_{s=0}^{\infty}M(s)$ be an admissible $V$-module. Then
$M(m)$ is an $A_m(V)$-module and $M(n)$ is an $A_n(V)$-module. As a result,
$\Hom_{\C}(M(m),M(n))$ is an $A_n(V)-A_m(V)$-bimodule. This observation
leads to the study of $A_n(V)-A_m(V)$-bimodules $A_{n,m}(V)$  in \cite{DJ}.
While the definition of $A_{n,m}(V)$ sounds complicated, the
representation-theoretical meaning of $A_{n,m}(V)$ is clear.

For homogeneous $u\in V,$ $v\in V$ and $m,n,p\in{\mathbb
Z}_{+}$, define the product $\ast_{m,p}^{n}$  on $V$ as follows
$$
u\ast_{m,p}^{n}v=\sum\limits_{i=0}^{p}(-1)^{i}{m+n-p+i\choose
i}{\rm Res}_{z}\frac{(1+z)^{wtu+m}}{z^{m+n-p+i+1}}Y(u,z)v.
$$
In order to explain the meaning of the product $u\ast_{m,p}^{n}v$ we consider
an admissible $V$-module $M=\sum_{s\geq 0}M(s).$ For $u\in V$
we set $o_t(u)=u_{\wt u-1-t}$ for $t\in \Z.$ Then $o_t(u)M(s)\subset M(s+t).$
It is proved in \cite{DJ} that $o_{n-p}(u)o_{p-m}(v)=o_{n-m}(u\ast_{m,p}^{n}v)$
acting on $M(m).$

If $n=p$, we denote $\ast_{m,p}^{n}$ by $\bar{\ast}_{m}^{n}$, and
if $m=p$, we denote $\ast_{m,p}^{n}$ by $\ast_{m}^{n}$, i.e.,
$$
u\ast_{m}^{n}v=\sum\limits_{i=0}^{m}(-1)^{i}{ n+i\choose i}{\rm
Res}_{z}\frac{(1+z)^{wtu+m}}{z^{n+i+1}}Y(u,z)v,
$$$$
u\bar{\ast}_{m}^{n}v=\sum\limits_{i=0}^{n}(-1)^{i}{ m+i\choose
i}{\rm Res}_{z}\frac{(1+z)^{wtu+m}}{z^{m+i+1}}Y(u,z)v.$$
From the discussion above we see that
 $o(u)o_{n-m}(v)=o_{n-m}(u\bar{\ast}_{m}^{n}v)$
and  $o_{n-m}(u)o(v)=o_{n-m}(u\ast_{m}^{n}v)$
on $M(m).$ So the products $u\ast_{m}^{n}v$ and $
u\bar{\ast}_{m}^{n}v$ will induce the right $A_m(V)$-module
and left $A_n(V)$-module structure on $A_{n,m}(V)$ which will
be defined later. It is clear that if $m=n,$ then  $u\ast_{m}^{n}v,$ $u\bar{\ast}_{m}^{n}v$
and  $u\ast_{n}v$ are equal.

Let $O'_{n,m}(V)$ be the linear span of all $u\circ_{m}^{n}v$ and
 $L(-1)u+(L(0)+m-n)u$, where for homogeneous $u\in V$ and $v\in V$,
$$
u\circ_{m}^{n}v={\rm
Res}_{z}\frac{(1+z)^{wtu+m}}{z^{n+m+2}}Y(u,z)v.$$ Again if $m=n,$
$u\circ_{m}^{n}v=u\circ_{n} v$ has been defined in \cite{DLM3}
(see Section 4). Then $O_{n}(V)=O'_{n,n}(V)$ (see \cite{DLM3} and
Section 4).

We introduce more
subspaces of $V.$ Let $O''_{n,m}(V)$ be the linear span of
$u\ast_{m,p_{3}}^{n}((a\ast_{p_{1},p_{2}}^{p_{3}}b){\ast}_{m,p_{1}}^{p_{3}}c-
a{\ast}_{m,p_{2}}^{p_{3}}(b{\ast}_{m,p_{1}}^{p_{2}}c)),$
 for $a,b,c,u\in V, p_{1},p_{2},p_{3}\in{\mathbb
Z}_{+}$, and $O'''_{n,m}(V)=\sum_{p\in{\mathbb
Z}_{+}}(V\ast_{p}^{n}O_{p}(V))\ast_{m,p}^{n}V.$ Set
$$
O_{n,m}(V)=O'_{n,m}(V)+O''_{n,m}(V)+O'''_{n,m}(V).
$$
and $$A_{n,m}(V)=V/O_{n,m}(V).$$

The following theorem about $A_{n,m}(V)$ is obtained in \cite{DJ}.

\begin{theorem}\label{t2.8} Let $V$ be a vertex operator algebra and $m,n$ nonnegative integers. Then

(1) $A_{n,m}(V)$ is an $A_n(V)$-$A_m(V)$-bimodule such that the
left and right actions of $A_n(V)$ and $A_m(V)$ are given by
$\bar*_m^n$ and $*^n_m.$

(2) $A_{m,n}(V)$ can be made to be an $A_{n}(V)-A_m(V)$-bimodule
isomorphic to $A_{n,m}(V).$

(3) Let $l$ be nonnegative integers such that $m-l,n-l$ are
nonnegative. Then $A_{n-l,m-l}(V)$ is an
$A_n(V)$-$A_m(V)$-bimodule and the identity map on $V$ induces an
epimorphism of $A_n(V)$-$A_m(V)$-bimodules from $A_{n,m}(V)$ to
$A_{n-l,m-l}(V).$

(4)  Define a linear map $\psi$: $A_{n,p}(V)\otimes_{A_p(V)}A_{p,m}(V)\rightarrow
A_{n,m}(V)$ by
$$
\psi(u\otimes v)=u\ast_{m,p}^{n}v,$$ for $u\otimes v\in
A_{n,p}(V)\otimes_{A_p(V)}A_{p,m}(V)$. Then $\psi$ is an
$A_{n}(V)-A_{m}(V)$- bimodule homomorphism from
$A_{n,p}(V)\otimes_{A_p(V)}A_{p,m}(V)$ to $A_{n,m}(V)$.

(5) Let $M=\sum_{s=0}^{\infty}M(s)$ be an admissible $V$-module.
Set $o_{n,m}(v)=v_{\wt(v)-1+m-n}.$ Then $v+O_{n,m}(V)\mapsto
o_{n,m}(v)$ gives an $A_n(V)-A_m(V)$-bimodule homomorphism from
$A_{n,m}(V)$ to $\Hom_{\C}(M(m),M(n)).$

(6) For any $n\geq 0,$ the $A_n(V)$ and $A_{n,n}(V)$ are the same.

(7) Let $U$ be an $A_{m}(V)$-module which can not factor through
$A_{m-1}(V)$. Then
$$
\bigoplus_{n\in{\mathbb Z}_{+}}A_{n,m}(V)\otimes_{A_{m}(V)}
U$$
is an admissible $V$-module isomorphic to $M(U)$ given in Theorem \ref{tha}
with $M(U)(n)=A_{n,m}(V)\otimes_{A_m(V)}U.$

(8) If $V$ is rational and $W^{j}\!=\!\bigoplus_{n\geq 0}\!W^j(n)$
with $W^j(0)\neq 0$ for $j=1,2,\cdots, s$ are all the inequivalent
irreducible modules of $V,$ then
$$A_{n,m}(V)\cong \bigoplus_{l=0}^{{\rm
min}\{m,n\}}\left(\bigoplus_{i=1}^{s}{\rm Hom}_{\mathbb
C}(W^i(m-l),W^i(n-l))\right).
$$
\end{theorem}

By the definitions $A_{n,n}(V)$ is a quotient of $A_n(V)$ as
$O_n(V)$ is a subspace of $O_{n,n}(V).$  The result (6) is  highly
nontrivial  and is proved by using the representation theory. It
is an interesting problem to find a direct proof  that $O_{n}(V)$
and $O_{n,n}(V)$ are equal.

So far we have two constructions of Verma type admissible module
$M(U)$ generated by an $A_m(V)$-module $U$ from Theorems \ref{tha}
and \ref{t2.8}. But the construction of $M(U)$ given in
Theorem \ref{tha} is a quotient module for certain Lie algebra
associated to $V$ (see \cite{DLM3} for the detail)
so the structure is not clear. On the other
hand, the construction of $M(U)$ from Theorem \ref{t2.8} (7) is
explicit and each homogeneous subspace $M(U)(n)$ is computable.
This new construction of $M(U)$ is fundamental in our study of
rationality in the next section.

\section{Rationality}
\setcounter{equation}{0}

In this section we discuss the relation between rationality of $V$
and semisimplicity of $A(V).$ This part is totally new and has not
appeared anywhere else.

Let $U$ be an $A(V)$-module. Then $M(U)=\oplus_{n\geq
0}A_{n,0}(V)\otimes_{A(V)}U$ is the Verma type admissible
$V$-module such that $M(U)(0)=U.$ It is easy to see that $M(U)$
has a maximal admissible submodule $J(U)=\sum_{n\geq 0}J(U)(n)$
such that $J(U)(0)=0$ where $J(U)(n)=M(U)(n)\cap J(U).$ As in
\cite{DLM2} we denote the quotient of $M(U)$ modulo $J(U)$ by
$L(U).$ Clearly, $L(U)$ is irreducible if and only if $U$ is
irreducible. From the representation theory for the affine
Kac-Moody algebras or the Virasoro algebra we know that in general
$M(U)$ and $L(U)$ are different. So it is natural to ask when
$L(U)$ and $M(U)$ are equal.

\begin{theorem}\label{t2.7} Let $V$
be a simple  vertex operator algebra such that $A(V)$ is
semisimple. Let $U$ be an irreducible module of $A(V)$, then the
Verma admissible $V$-module $M(U)=\bigoplus_{n\in{\mathbb
Z}_{+}}A_{n,0}(V)\otimes U$ generated by $U$ is an irreducible
admissible module of $V$. That is, $M(U)=L(U).$
\end{theorem}

The proof of Theorem \ref{t2.7} is highly nontrivial and will be given
in another paper. The main idea in the proof is to use the bimodule
theory developed in \cite{DJ} to show that
$M(U)(n)=A_{n,0}(V)\otimes U$ is an irreducible $A_n(V)$-module for
all $n.$

It is proved in \cite{DLTYY} that if $M(U)=L(U)$ for any simple
$A(V)$-module, $A(V)$ is semisimple and $V$ is $C_2$-cofinite then
$V$ is rational. The same proof works here with $C_2$-cofiniteness
replaced by the assumption that each irreducible admissible
$V$-module is ordinary. So using Theorem \ref{t2.7} we have the
following theorem.
\begin{theorem}\label{tt2.8} Let $V$
be a simple  vertex operator algebra. Then $V$ is rational if and only if
 $A(V)$ is semisimple and each irreducible admissible $V$-module is ordinary.
\end{theorem}

The assumption that any irreducible admissible module is ordinary
in Theorem \ref{tt2.8} is equivalent to that each homogeneous
subspace of any irreducible admissible module is finite
dimensional. It is worthy to point out that this is not a very
strong assumption which is true for all known simple vertex
operator algebras.

\begin{coro} Let $V$ be a $C_2$-cofinite simple vertex operator algebra.
Then the following are equivalent:

(1) $V$ is regular,

(2) $V$ is rational,

(3) $A(V)$ is semisimple.
\end{coro}

The corollary follows immediately from the fact that any finitely generated
admissible module is ordinary (see \cite{KL}, \cite{Bu}) if $V$ is
$C_2$-cofinite. We have already mentioned that rationality together with
the $C_2$-cofiniteness and regularity are equivalent.  In the proof
of modular invariance of trace functions in \cite{Z}, the semisimplicity
of $A(V)$ (not the rationality of $V$) and $C_2$-cofiniteness were used.
But under the assumption of $C_2$-cofiniteness, the rationality of $V$
and semisimplicity of $A(V)$ are equivalent. So the rationality of $V$ is,
in fact, used in the proof of modularity in \cite{Z}.

Recall from \cite{KL} that $V$ is called $C_1$-cofinite if
$V=\sum_{n\geq 0}V_n$ with $V_0=\C \1$ such that $C_1(V)$ has
finite codimension in $V$ where $C_1(V)$ is spanned by $u_{-1}v$
and $L(-1)u$ for all $u,v\in \sum_{n>0}V_n.$ It is proved in
\cite{KL} that if $V$ is $C_1$-cofinite then any irreducible
admissible $V$-module is ordinary. So we have another important
corollary:

\begin{coro} Let $V$ be a $C_1$-cofinite simple vertex operator algebra.
Then $V$ is rational if and only if $A(V)$ is semisimple.
\end{coro}

Now we can easily conclude the rationality of many vertex operator
algebras which were very hard theorems. Here is an example.
Let $L$ be a positive
definite even lattice and $V_L$ the corresponding vertex operator
algebra (see \cite{B} and \cite{FLM}). Then $V_L$ has a canonical
automorphism $\theta$ of order 2 induced from the $-1$ isometry of
the lattice and the fixed point subspace $V_L^+$ is a simple
vertex operator algebra. The irreducible admissible modules for $V_L^+$
have been classified in \cite{DN2} and \cite{AD}. Moreover,
each irreducible admissible module is ordinary (see \cite{DN2}, \cite{AD}).
\begin{coro} $A(V_L^+)$ is semisimple and $V_L^+$ is rational for any positive
definite even lattice $L.$
\end{coro}

We remark that if the rank of $L$ is one, the semisimplicity of
$A(V_L^+)$ has been proved in \cite{DN2} and the rationality of
$V_L^+$  has been established in \cite{A1} using a different
method.

We end this paper with the following conjecture.
\begin{conjecture} Let $V$ be a simple vertex operator algebra. Then
$V$ is rational if and only if $A(V)$ is semisimple.
\end{conjecture}

\end{document}